\documentclass{amsart}
\usepackage{amssymb,latexsym}
\input epsf
\newtheorem{theorem}{Theorem}[section]
\newtheorem{proposition}[theorem]{Proposition}
\newtheorem{corollary}[theorem]{Corollary}

\newtheorem{remark}[theorem]{Remark}
\newtheorem{lemma}[theorem]{Lemma}
\newtheorem{example}[theorem]{Example}

\newcommand{\lk}{{\rm lk}}
\newcommand{\nc}{{\rm {\bf L}}}
\newcommand{\KK}{{\mathbb K}}
\newcommand{\RR}{{\mathbb R}}

\newcommand{\ZZ}{{\mathbb Z}}

\newcommand{\iI}{{\mathcal I}}

\newcommand{\rR}{{\mathcal R}}

\newcommand{\ppP}{{\bf P}}

\renewcommand{\to}{\rightarrow}

\newcommand{\sm}{{\setminus}}
\begin{document}
\title[Shellability of generalized cluster complexes]
{Shellability and higher Cohen-Macaulay connectivity of generalized
cluster complexes}

\author{Christos~A.~Athanasiadis}
\address{Department of Mathematics (Division of Algebra-Geometry)\\
University of Athens\\
Panepistimioupolis\\
15784 Athens, Greece}
\email{caath@math.uoa.gr}

\author{Eleni~Tzanaki}
\address{Department of Mathematics\\
University of Crete\\
71409 Heraklion, Crete, Greece}
\email{etzanaki@math.uoc.gr}
\date{July 3, 2006; revised, February 27, 2007}
\begin{abstract}
Let $\Phi$ be a finite root system of rank $n$ and let $m$ be a
nonnegative integer. The generalized cluster complex $\Delta^m (\Phi)$ 
was introduced by S. Fomin and N. Reading. It was conjectured by 
these authors that $\Delta^m (\Phi)$ is shellable and by V. Reiner 
that it is $(m+1)$-Cohen-Macaulay, in the sense of Baclawski. 
These statements are proved in this paper. Analogous statements are 
shown to hold for the positive part $\Delta^m_+ (\Phi)$ of $\Delta^m 
(\Phi)$. An explicit homotopy equivalence is given between $\Delta^m_+ 
(\Phi)$ and the poset of generalized noncrossing partitions, associated 
to the pair $(\Phi, m)$ by D. Armstrong.
\end{abstract}

\maketitle

\section{Introduction}
\label{intro}

Let $\Phi$ be a finite root system of rank $n$ and let $m$ be a 
nonnegative integer. The generalized cluster complex $\Delta^m (\Phi)$ 
was introduced by S. Fomin and N. Reading \cite{FR2}, soon after it was
studied in \cite{Tz1} in the special cases in which $\Phi$ has type $A$ 
or $B$ in the Cartan-Killing classification. If $\Phi$ has type $A$, then
faces of $\Delta^m (\Phi)$ can be described in terms of polygonal 
subdivisions of a convex polygon with $m(n+1)+2$ vertices (see Section 
\ref{pre2}). If $m=1$ and $\Phi$ has type $A$, then $\Delta^m 
(\Phi)$ is combinatorially isomorphic to the boundary complex of the 
$n$-dimensional simplicial 
associahedron, a classical object of study in enumerative and polyhedral 
combinatorics. More generally, $\Delta^m (\Phi)$ reduces to the cluster 
complex $\Delta (\Phi)$ when $m=1$, a simplicial complex of importance 
in the context of cluster algebras and $Y$-systems \cite{FZ1, FZ2}; 
see \cite{FR1} for a nice introduction to these topics.

Motivation for defining and studying generalized cluster complexes came
from at least two directions. In \cite{FR2} combinatorial algorithms for
determining Coxeter theoretic invariants were given, in which certain
identities satisfied by the face numbers of generalized cluster complexes
are crucial. On the other hand, it was conjectured by the first author that
the $m$-generalized Narayana (or Fuss-Narayana) numbers defined in
\cite{Ath1} for a crystallographic root system $\Phi$ form the $h$-vector
of an $(n-1)$-dimensional Cohen-Macaulay simplicial complex. It follows 
from the results of this paper (as well as those of \cite{Ath1, FR2}) 
that the generalized cluster complex $\Delta^m (\Phi)$ is the desired 
complex having these properties (this was verified earlier by the second 
author in the special cases of \cite{Tz1}). Moreover, it has been shown
\cite{Ar2, AT, FR2, Tz1, Tz2} that the complex $\Delta^m (\Phi)$, as
well as a natural subcomplex $\Delta^m_+ (\Phi)$ called its positive
part, share essentially all enumerative properties of cluster complexes,
relating to a variety of interesting structures within the algebraic
combinatorics of Coxeter groups, root systems and hyperplane arrangements;
see \cite{Ath1} \cite[Lecture 5]{FR1} \cite{FR2}.

This paper shows that the complexes
$\Delta^m (\Phi)$ and $\Delta^m_+ (\Phi)$ have attractive topological
properties as well. Recall that a simplicial complex $\Delta$ is
said to be \emph{$k$-Cohen-Macaulay} (over $\ZZ$ or some field $\KK$) if
the complex obtained from $\Delta$ by removing any subset of its vertex
set of cardinality less than $k$ is Cohen-Macaulay (over $\ZZ$ or $\KK$)
and has the same dimension as $\Delta$. The following theorem is the main
result of this paper.

\begin{theorem}
\begin{enumerate}
\itemsep=0pt \item[(i)] The simplicial complex obtained from
$\Delta^m (\Phi)$ by removing any subset of its vertex set of
cardinality not exceeding $m$ is pure, of dimension $n-1$, and
shellable. In particular, $\Delta^m (\Phi)$ is shellable, hence
homotopy equivalent to a wedge of $(n-1)$-dimensional spheres, and
$(m+1)$-Cohen-Macaulay.

\item[(ii)] The simplicial complex obtained from $\Delta^m_+ (\Phi)$
by removing any subset of its vertex set of cardinality not exceeding
$m-1$ is pure, of dimension $n-1$, and shellable. In particular,
$\Delta^m_+ (\Phi)$ is pure of dimension $n-1$ and shellable, hence
homotopy equivalent to a wedge of $(n-1)$-dimensional spheres, and
$m$-Cohen-Macaulay.
\end{enumerate}
\label{thm0}
\end{theorem}

A few comments on the theorem are in order. Shellability of $\Delta^m
(\Phi)$ was conjectured by S. Fomin and N. Reading \cite[Conjecture
11.3]{FR2} (see also \cite[Problem 6.3]{Ar1}) and verified by the second
author \cite{Tz1} when $\Phi$ has type $A$ or $B$. This conjecture was
extended to $\Delta^m_+ (\Phi)$ by the authors \cite[Conjecture
4.6]{AT}. The question of higher Cohen-Macaulay connectivity of $\Delta^m
(\Phi)$ was raised by V. Reiner; see \cite[p. 17]{Ar1} (later, the same
question was raised independently by G. Kalai \cite{Ka}).
The result of Theorem \ref{thm0} in this respect is best possible (see
Remark \ref{rem:best}) and nontrivial even when $\Phi$ has type $A$. The 
Euler characteristic of $\Delta^m (\Phi)$ was computed in \cite[Section 
11]{FR2}.
It is known \cite[Section 10]{FR2} that the $m$-generalized Narayana
numbers defined in \cite{Ath1} for crystallographic $\Phi$ coincide with
the entries of the $h$-vector of $\Delta^m (\Phi)$, except possibly when
$m \ge 2$ and $\Phi$ contains an irreducible component of type $F_4$,
$E_6$, $E_7$ or $E_8$. Hence the statement on shellability in Theorem
\ref{thm0} (i) establishes the conjecture on these numbers, mentioned
earlier, in these cases. 

The concept of higher Cohen-Macaulay connectivity was introduced by K.
Baclawski \cite{Ba}. We refer the reader to \cite{Ba} \cite{Sw}
\cite[Section III.3]{St}
and references given in these sources for other interesting classes of
simplicial complexes known to be $k$-Cohen-Macaulay for some $k \ge 2$.

The general layout of this paper is as follows. Section \ref{pre}
includes background on simplicial complexes, (generalized) cluster
complexes and a related partial order on a finite real reflection
group. In particular, a new characterization (Theorem \ref{thm:Tz}) of
$\Delta^m (\Phi)$, due to the second author \cite{Tz2, Tz3}, is reviewed,
generalizing the one for $\Delta (\Phi)$ given by T. Brady and C. Watt
\cite[Section 8]{BW}. The proof of Theorem \ref{thm0}, which relies on
this characterization, is given in Section \ref{proof}. That section
includes also a computation of the Euler characteristic of $\Delta^m_+
(\Phi)$ (Corollary \ref{cor:euler}). Finally some applications to the
topology of the posets of generalized noncrossing partitions, introduced
and studied by D. Armstrong \cite{Ar2}, are given in Section \ref{gnc}.

\section{Background}
\label{pre}

Throughout the paper we denote by $|\sigma|$ the cardinality of a
finite set $\sigma$.

\subsection{Simplicial complexes.}
\label{pre1}
Let $E$ be a finite set. An (abstract) \emph{simplicial complex} on the
ground set $E$ is a collection $\Delta$ of subsets of $E$ such that $\tau
\subseteq \sigma \in \Delta$ implies $\tau \in \Delta$. The set $V = \{v
\in E: \{v\} \in \Delta\}$ is the set of \emph{vertices} of $\Delta$. The
elements of $\Delta$ are called \emph{faces} and those which are maximal
with respect
to inclusion are called \emph{facets}. The dimension of a face $\sigma$ is
defined as one less than the cardinality of $\sigma$ and the dimension of
$\Delta$ as the maximum dimension of a face. The complex $\Delta$ is
\emph{pure} if all its facets have the same dimension and \emph{flag} if
all its minimal nonfaces have two elements. The
$k$-skeleton $\Delta^{\le k}$ of $\Delta$ is the subcomplex
formed by the faces of $\Delta$ of dimension at most $k$. The simplicial
join $\Delta_1 \ast \Delta_2$ of two abstract simplicial complexes
$\Delta_1$ and $\Delta_2$ on disjoint ground sets has as its faces the 
sets of the form $\sigma_1 \cup \sigma_2$, where $\sigma_1 \in \Delta_1$ 
and $\sigma_2 \in \Delta_2$.

The \emph{link} of $v \in E$ in $\Delta$ is defined as $\lk_\Delta
(v) = \{\sigma \sm \, \{v\}: \sigma \in \Delta, \, v \in \sigma\}$. The 
\emph{induced subcomplex} or \emph{restriction} of $\Delta$ on $A 
\subseteq E$ is defined as $\Delta_A = \{\sigma \in \Delta: \sigma 
\subseteq A\}$. We will write $\Delta \sm v$ for the restriction of 
$\Delta$ on the set $E \sm \, \{v\}$. A pure simplicial complex 
$\Delta$ is \emph{shellable} if there exists a total ordering $\sigma_1, 
\sigma_2,\dots,\sigma_m$ of the set of facets of $\Delta$ such that 
for any given indices $1 \le i < k \le m$ there exists $1 \le j < k$
and $v \in \sigma_k$ with $\sigma_i \cap \sigma_k \subseteq \sigma_j 
\cap \sigma_k = \sigma_k \sm \, \{v\}$. The following lemma is 
elementary and well known; see, for instance, \cite[Section 2]{PB}.
\begin{lemma}
\begin{enumerate}
\itemsep=0pt
\item[(i)]
If $\Delta \sm v$ is pure of dimension $d$ and shellable and $\lk_\Delta
(v)$ is pure of dimension $d-1$ and shellable for some $v \in E$ then
$\Delta$ is pure of dimension $d$ and shellable.

\item[(ii)]
The simplicial join of pure shellable complexes is pure and shellable.
\qed
\end{enumerate}
\label{lem0}
\end{lemma}
When we talk about algebraic or topological properties of an abstract
simplicial complex $\Delta$ we refer to those of its Stanley-Reisner 
ring \cite[Chapter II]{St} or geometric realization \cite[Section 
9]{Bj}, respectively. Thus any pure $d$-dimensional, shellable 
simplicial complex $\Delta$ is Cohen-Macaulay over $\ZZ$ and all 
fields and homotopy equivalent to a wedge of spheres of dimension 
$d$.

\subsection{Generalized cluster complexes.}
\label{pre2}
Let $\Phi$ be a finite root system spanning $\RR^n$, endowed with
the standard inner product, and let $\Phi^+$ be a fixed choice of
a positive system with corresponding simple system $\Pi$. We assume 
temporarily that $\Phi$ is irreducible and let $\Pi = 
\{\alpha_1,\dots,\alpha_n\}$ be ordered so that for some $1 \le s 
\le n$, each of the sets $\Pi_+ = \{\alpha_1,\dots,\alpha_s\}$ and 
$\Pi_- = \{\alpha_{s+1},\dots,\alpha_n\}$ has pairwise orthogonal 
elements. Let $\Phi_{\ge -1} = \Phi^+ \cup (-\Pi)$ and define the 
involutions $\tau_\pm: \Phi_{\ge -1} \to \Phi_{\ge -1}$ by   
\[ \tau_\varepsilon (\alpha) \ = \
\begin{cases} \alpha, & \text{if $\alpha \in -\Pi_{-\varepsilon}$} \\
R_\varepsilon (\alpha), & \text{otherwise} \end{cases} \]
for $\varepsilon \in \{+,-\}$ and $\alpha \in \Phi_{\ge -1}$, where 
$R_\varepsilon$ is the product of the reflections (taken in any order)
in the linear hyperplanes in $\RR^n$ orthogonal to the elements of 
$\Pi_\varepsilon$. The product $\rR = \tau_- \tau_+$, introduced in 
\cite{FZ2}, can be viewed as a deformation of the Coxeter element in 
the real reflection group corresponding to $\Phi$. Given a nonnegative 
integer $m$, the set
\[ \Phi^m_{\ge -1} = \{\alpha^i: \alpha \in \Phi^+ \
\text{and} \ 1 \le i \le m\} \, \cup \, (-\Pi) \]
consists of the negative simple roots and $m$ copies of each positive
root, each copy colored with one of $m$ possible colors. Using the
convention $\alpha^1 = \alpha$ for $\alpha \in -\Pi$, the map $\rR_m: 
\Phi^m_{\ge -1} \to \Phi^m_{\ge -1}$ is defined by
\[ \rR_m (\alpha^k) \ = \
\begin{cases} \alpha^{k+1}, 
& \text{if $\alpha \in \Phi^+$ \text{and} $k<m$} \\
(\rR (\alpha))^1, & \text{otherwise.} \end{cases} \]
It can be shown \cite[Theorem 3.4]{FR2} that there exists a unique
symmetric binary relation on $\Phi^m_{\ge -1}$, called ``compatibility",
which has the following two properties:  
\begin{enumerate}
\itemsep=0pt
\item[$\bullet$]
$\alpha^k$ is compatible with $\beta^\ell$ if and only if $\rR_m 
(\alpha^k)$ is compatible with $\rR_m (\beta^\ell)$; 
\item[$\bullet$]
for $\alpha \in \Pi$, the root $-\alpha$ is compatible with $\beta^\ell$ 
if and only if the simple root expansion of $\beta$ does not involve 
$\alpha$.
\end{enumerate}
The generalized cluster complex $\Delta^m (\Phi)$ is defined in 
\cite[Section 3]{FR2} is the abstract simplicial complex on the vertex 
set $\Phi^m_{\ge -1}$ which has as its faces the subsets of mutually
compatible elements of $\Phi^m_{\ge -1}$. If $\Phi$ is a direct product 
$\Phi_1 \times \Phi_2$ then $\Delta^m (\Phi)$ is defined as the simplicial
join of $\Delta^m (\Phi_1)$ and $\Delta^m (\Phi_2)$. 

The complex $\Delta^m (\Phi)$ is flag by definition and pure of 
dimension $n-1$ \cite[Theorem 3.9]{FR2}. 
Following \cite{AT}, we denote by $\Delta^m_+ (\Phi)$ the induced
subcomplex of $\Delta^m (\Phi)$ on the set of vertices obtained from
$\Phi^m_{\ge -1}$ by removing the negative simple roots and call this
simplicial complex the \emph{positive part} of $\Delta^m (\Phi)$. For
$\alpha \in \Pi$ we denote by $\Phi_\alpha$ the standard parabolic
root subsystem obtained by intersecting $\Phi$ with the linear span of
$\Pi \sm \, \{\alpha\}$, endowed with the induced positive system
$\Phi_\alpha^+ = \Phi^+ \cap \Phi_\alpha$. Let us summarize the 
properties of $\Delta^m (\Phi)$ which will be of importance for us.
\begin{proposition} {\rm (\cite{FR2})}
\begin{enumerate}
\itemsep=0pt
\item[(i)] If $\Phi$ is a direct product $\Phi_1 \times \Phi_2$
then $\Delta^m (\Phi) = \Delta^m (\Phi_1) \ast \Delta^m (\Phi_2)$ and
$\Delta^m_+ (\Phi) = \Delta^m_+ (\Phi_1) \ast \Delta^m_+ (\Phi_2)$.
\end{enumerate}
Suppose that $\Phi$ is irreducible.
\begin{enumerate}
\itemsep=0pt
\item[(ii)] If $\alpha \in \Pi$ then $\lk_{\Delta^m (\Phi)} (-\alpha) =
\Delta^m (\Phi_\alpha)$.

\item[(iii)] The map $\rR_m: \Phi^m_{\ge -1} \to \Phi^m_{\ge -1}$ is 
a bijection with the following properties:
\begin{itemize}
\itemsep=0pt
\item[(a)] $\sigma \in \Delta^m (\Phi)$ if and only if $\rR_m (\sigma) :=
\{\rR_m (\alpha): \alpha \in \sigma\} \in \Delta^m (\Phi)$.
\item[(b)] For any $\beta \in \Phi^m_{\ge -1}$ there exists $j$ such
that $\rR_m^j (\beta) \in (-\Pi)$. \qed
\end{itemize}
\end{enumerate}
\label{prop0}
\end{proposition}

It is not hard to describe $\Delta^m (\Phi)$ explicitly when $\Phi$ has
type $A_{n-1}$ \cite{FR2, Tz1}. Call a diagonal of a convex polygon $\ppP$
with $mn+2$ vertices \emph{$m$-allowable} if it divides $\ppP$ into
two polygons, each with number of vertices congruent to $2 \mod m$.
Then vertices of $\Delta^m (\Phi)$ biject to the $m$-allowable
diagonals of $\ppP$ so that faces correspond to the sets of pairwise
noncrossing diagonals of this kind.

\subsection{The reflection length order.}
\label{pre3}
We will denote by $R(\alpha)$ the reflection in the linear hyperplane
in $\RR^n$ orthogonal to a nonzero (colored or not) vector $\alpha \in
\RR^n$. Let $W$ be the finite real reflection group generated by the
set $T$ of reflections $R(\alpha)$ for $\alpha \in \Phi$. For $w \in
W$ we denote by $\ell_T (w)$ the smallest integer $r$ such that
$w$ can be written as a product of $r$ reflections in $T$. The set $W$
can be partially ordered by letting
\[ u \preceq w \ \ \ \text{if and only if} \ \ \ \ell_T (u) + \ell_T
(u^{-1} w) = \ell_T (w), \]
in other words if there exists a shortest factorization of $u$ into
reflections in $T$ which is a prefix of such a shortest factorization
of $v$. The order $\preceq$, called \emph{reflection length order}
or \emph{absolute order}, turns $W$ into a graded poset having the
identity of $W$ as its unique minimal element and rank function
$\ell_T$. If $\Phi$ has type $A_{n-1}$, so that $W$ can be realized 
as the symmetric group of permutations of the set $\{1, 2,\dots,n\}$,
then one can describe more explicitly these concepts as follows. We 
have $\ell_T (w) = n - c(w)$, where $c(w)$ is the number of cycles 
of $w$, and $u \preceq w$ if and only if each cycle of $u$ can be 
obtained from some cycle of $w$ by deleting elements.  

We will assume elementary properties of the order $\preceq$;
see \cite[Section 2]{BW} \cite[Section 2]{ABMW}. For instance, if
$v \preceq u \preceq w$ then $v^{-1} u \preceq v^{-1} w$ and $u
v^{-1} \preceq w v^{-1}$.

\subsection{Generalized cluster complexes via the reflection length
order.} \label{pre4}
Suppose that $\Phi$ is irreducible and let $N$ denote the cardinality
of $\Phi^+$. As in Section \ref{pre2}, let $\Pi = 
\{\alpha_1,\dots,\alpha_n\}$ be ordered so that for some $1 \le s \le 
n$ the sets $\Pi_+ = \{\alpha_1,\dots,\alpha_s\}$ and $\Pi_- = 
\{\alpha_{s+1},\dots,\alpha_n\}$ have pairwise orthogonal elements 
and let
\[ \gamma = R(\alpha_1) R(\alpha_2) \cdots R(\alpha_n) \]
be a corresponding bipartite Coxeter element of $W$. As in \cite{ABMW,
BW} we set $\rho_i = R(\alpha_1) R(\alpha_2) \cdots R(\alpha_{i-1})
(\alpha_i)$ for $i \ge 1$, where the $\alpha_i$ are indexed cyclically
modulo $n$ (so that $\rho_1 = \alpha_1$), and $\rho_{-i} = \rho_{2N-i}$
for $i \ge 0$ and recall that
\[ \begin{tabular}{ll}
$\{ \rho_1, \rho_2,\dots,\rho_N \} = \Phi^+$, \\
$\{ \rho_{N+i}: 1 \le i \le s\} = \{-\rho_1,\dots,-\rho_s \} =
- \Pi_+$, \\
$ \{\rho_{-i}: 0 \le i < n-s \} = \{-\rho_{N-i}: 0 \le i < n-s\} =
- \Pi_-$.
\end{tabular} \]
We consider the total order $<$ on the set $\Phi_{\ge -1} = \Phi^+ 
\cup (-\Pi)$ defined by
\begin{equation}
\rho_{-n+s+1} < \cdots < \rho_0 < \rho_1 < \cdots < \rho_{N+s}.
\label{total}
\end{equation}
This order induces a total order, which we denote again by $<$, on the
set of elements of $\Phi^m_{\ge -1}$ which are positive roots of some
fixed color $i$ simply by forgetting the color. For $\tau = \{\tau_1,
\tau_2,\dots,\tau_k\} \subseteq \Phi^m_{\ge -1}$ such that either
$\tau \subseteq (-\Pi)$ or $\tau$ consists of positive roots of the
same color, we let
\begin{equation}
w_\tau = R(\tau_k) R(\tau_{k-1}) \cdots R(\tau_1),
\label{eq:w_tau}
\end{equation}
where $\tau_1 < \tau_2 < \cdots < \tau_k$, with the convention that
$w_\tau$ is the identity element of $W$ if $\tau = \emptyset$. For any
$\sigma \subseteq \Phi^m_{\ge -1}$ we let
\begin{equation}
w_\sigma = w_{\sigma^+} w_{\sigma^{(m)}} w_{\sigma^{(m-1)}} \cdots
w_{\sigma^{(1)}} w_{\sigma^-}
\label{eq:w_sigma}
\end{equation}
where $\sigma^\pm = (-\Pi_\pm) \cap \sigma$ and for $1 \le i \le m$,
$\sigma^{(i)}$ is the set of elements of $\sigma$ which are positive
roots of color $i$. The following characterization of $\Delta^m (\Phi)$
generalizes that of $\Delta (\Phi)$ given in \cite[Section 8]{BW}.
\begin{theorem} {\rm (\cite{Tz2} \cite[Section 4.3]{Tz3})}
A set $\sigma \subseteq \Phi^m_{\ge -1}$ is a face of $\Delta^m (\Phi)$
if and only if $w_\sigma \preceq \gamma$ and $\ell_T (w_\sigma) =
|\sigma|$. \qed
\label{thm:Tz}
\end{theorem}
\begin{example} {\rm Suppose that $\Phi$ has type $A_2$. Then the simple 
roots $\alpha_1, \alpha_2$ can be chosen as unit vectors in the Euclidean 
plane forming an angle of $2 \pi / 3$ and we have $\Pi_+ = \{\alpha_1\}$, 
$\Pi_- = \{\alpha_2\}$ and $\Phi^+ = \{\alpha_1, \alpha_2, \alpha\}$, 
where $\alpha = \alpha_1 + \alpha_2$. Moreover, $W$ is the dihedral group 
of order 6 generated by the reflections $R(\alpha_1)$, $R(\alpha_2)$ and 
$R(\alpha)$ and $\gamma = R(\alpha_1) R(\alpha_2)$ is a rotation by $2 
\pi / 3$. One can compute easily that $\rho_1 = \alpha_1$, $\rho_2 = 
\alpha$, $\rho_3 = \alpha_2$, $\rho_4 = -\alpha_1$ and $\rho_0 = \rho_6
= -\alpha_2$. Therefore the total order (\ref{total}) on the set 
$\Phi_{\ge -1}$ is given by
\[ -\alpha_2 < \alpha_1 < \alpha < \alpha_2 < -\alpha_1. \]
Assume that $m=1$ and let $\tau = \{\alpha, \alpha_2\}$. We have $w_\tau 
=  R(\alpha_2) R(\alpha) = \gamma$ and hence, according to the condition 
in Theorem \ref{thm:Tz}, $\tau$ must be a facet of $\Delta (\Phi)$.
Similarly we find that $\Delta (\Phi)$ must have exactly four more facets,
namely $\{\alpha_1, \alpha\}$, $\{-\alpha_2, \alpha_1\}$, $\{-\alpha_1, -
\alpha_2\}$ and $\{\alpha_2, -\alpha_1\}$, corresponding to the 
factorizations $$R(\alpha) R(\alpha_1) = R(\alpha_1) R(-\alpha_2) = 
R(-\alpha_1) R(-\alpha_2) = R(-\alpha_1) R(\alpha_2)$$ of $\gamma$.
Indeed, $\Delta (\Phi)$ is the one-dimensional simplicial complex having
these five two-element sets as its facets \cite[p.~985]{FZ2}. Assume 
now that $m=2$ and $\tau = \{\alpha_2^1, \alpha_1^2\}$. We have $w_\tau 
= R(\alpha_1^2) R(\alpha_2^1) = R(\alpha_1) R(\alpha_2) = \gamma$ and 
hence, by Theorem \ref{thm:Tz}, $\tau$ is a facet of $\Delta^2 (\Phi)$.
Similarly we find that $\Delta^2 (\Phi)$ has the following 12 two-element 
sets as its facets: 
$\{\alpha_1^1, \alpha^1\}$, 
$\{\alpha_1^1, \alpha^2\}$, 
$\{\alpha_1^2, \alpha^2\}$, 
$\{\alpha^1, \alpha_2^1\}$,
$\{\alpha^1, \alpha_2^2\}$,
$\{\alpha^2, \alpha_2^2\}$,
$\{\alpha_2^1, \alpha_1^2\}$,  
$\{-\alpha_2, \alpha_1^1\}$, 
$\{-\alpha_2, \alpha_1^2\}$, 
$\{\alpha_2^1, -\alpha_1\}$, 
$\{\alpha_2^2, -\alpha_1\}$ and $\{-\alpha_1, - \alpha_2\}$. 
\qed}
\label{ex:A2}
\end{example}
Given $w \preceq \gamma$, we denote by $\Delta^m_+ (w)$ the subcomplex of
$\Delta^m_+ (\Phi)$ consisting of those faces $\sigma$ for which $w_\sigma
\preceq w$. Clearly the dimension of $\Delta^m_+ (w)$ is at most $\ell_T (w)
- 1$. It follows from \cite[Lemma 2.1 (iv)]{ABMW} that $\Delta^m_+ (w)$
coincides with the induced subcomplex of $\Delta^m (\Phi)$ on the vertex
set of positive colored roots $\alpha \in \Phi^m_{\ge -1}$ satisfying
$R(\alpha) \preceq w$. We will denote this vertex set by $\Phi^m_+ (w)$, so
that $\Phi^m_+ (\gamma)$ is the vertex set of $\Delta^m_+ (\Phi)$. When
$m=1$ the complex $\Delta^m_+ (w)$, written simply as $\Delta_+ (w)$, is
homeomorphic to a triangulation of a simplex of dimension $\ell_T (w)-1$
\cite[Corollary 7.7]{BW} (in particular it is pure of that dimension) and
was shown to be shellable in \cite[Theorem 7.1]{ABMW}.
\begin{theorem} {\rm (\cite{ABMW, BW})}
For any $w \preceq \gamma$ the complex $\Delta_+ (w)$ is pure of
dimension $\ell_T (w)-1$ and shellable. Moreover, it is homeomorphic to
a ball.
\qed
\label{thm:ABMW}
\end{theorem}

Our definitions of $\Delta^m_+ (w)$ and $\Phi^m_+ (w)$ extend naturally
when $\Phi$ is reducible. More precisely, if $\Phi = \Phi_1 \times \Phi_2$
then $\gamma = \gamma_1 \gamma_2$, where $\gamma_i$ is a bipartite Coxeter
element for the reflection group $W_i$ corresponding to $\Phi_i$, and $w
\preceq \gamma$ if and only if $w = w_1 w_2$ with $w_1 \preceq \gamma_1$
and $w_2 \preceq \gamma_2$. For $\sigma \in \Delta^m_+ (\Phi)$ we can
write $\sigma = \sigma_1 \cup \sigma_2$ with $\sigma_i \in \Delta^m_+
(\Phi_i)$ for $i = 1, 2$ and define $w_\sigma = w_{\sigma_1} w_{\sigma_2}$.
The subcomplex $\Delta^m_+ (w)$ of $\Delta^m_+ (\Phi)$ consisting of
those faces $\sigma$ with $w_\sigma \preceq w$ is equal to the simplicial
join of $\Delta^m_+ (w_1)$ and $\Delta^m_+ (w_2)$ and has vertex set
$\Phi^m_+ (w) = \Phi^m_+ (w_1) \cup \Phi^m_+ (w_2)$.

\section{Proof of Theorem \ref{thm0}}
\label{proof}

Theorem \ref{thm0} will be derived from the following proposition.
\begin{proposition}
Fix an index $1 \le j \le m$.
\begin{enumerate}
\itemsep=0pt
\item[(i)] Given $w \preceq \gamma$, the induced
subcomplex of $\Delta^m_+ (w)$ on any subset of its vertex set $\Phi^m_+
(w)$ containing all positive roots in $\Phi^m_+ (w)$ of color $j$ is pure
of dimension $\ell_T (w) - 1$ and shellable.

\item[(ii)] The induced subcomplex of $\Delta^m (\Phi)$ on any
subset of the vertex set of $\Delta^m (\Phi)$ containing all positive
roots of color $j$ is pure of dimension $n-1$ and shellable.
\end{enumerate}
\label{prop:pureshell}
\end{proposition}
\begin{proof}
(i) In view of the discussion following Theorem \ref{thm:ABMW} and
Lemma \ref{lem0} (ii), we may assume that $\Phi$ is irreducible.
Let $r = \ell_T (w)$ and recall that the dimension of $\Delta^m_+ (w)$
does not exceed $r-1$. Let $\Delta^m_A (w)$ denote the induced subcomplex
of $\Delta^m_+ (w)$ on the vertex set $A \subseteq \Phi^m_+ (w)$. Suppose
that $A$ contains all positive roots in $\Phi^m_+ (w)$ of color $j$. To
simplify notation we set $\Delta = \Delta^m_A (w)$ and let $\tau$ be any
face of $\Delta$. By Theorem \ref{thm:Tz} we have
\[ w_\tau = w_{\tau^{(m)}} w_{\tau^{(m-1)}} \cdots w_{\tau^{(1)}}
\preceq w \]
in the notation of (\ref{eq:w_sigma}), and $\ell_T (w_\tau) = |\tau|$.
It follows that $w_{\tau^{(j)}} \preceq u$ where
\[ u = (w_{\tau^{(m)}} \cdots w_{\tau^{(j+1)}})^{-1} \, w \,
(w_{\tau^{(j-1)}} \cdots w_{\tau^{(1)}})^{-1} \preceq w \]
and $\ell_T (u) = r - |\tau| + |\tau^{(j)}|$. Observe that the set
obtained from $\tau^{(j)}$ by uncoloring its elements is a face of
$\Delta_+ (u)$, which is pure of dimension $\ell_T (u) - 1$ by
Theorem \ref{thm:ABMW}. Hence there exists a set $\sigma^{(j)}$
containing $\tau^{(j)}$ and consisting of $\ell_T (u)$ positive
roots of color $j$ such that $w_{\sigma^{(j)}} = u$. Setting
$\sigma^{(i)} = \tau^{(i)}$ for all $1 \le i \le m$ with $i \neq
j$, we obtain a subset $\sigma$ of the vertex set of $\Delta^m_+
(\Phi)$ such that $\tau \subseteq \sigma$, $w_\sigma = w$ and
$|\sigma| = \ell_T (w) = r$. It follows from Theorem \ref{thm:Tz}
that $\sigma$ is a face of $\Delta^m (\Phi)$, and hence of
$\Delta^m_+ (w)$. Notice that (i) $A$ contains all positive roots 
in $\Phi^m_+ (w)$ of color $j$, (ii) $\sigma \sm \, \tau =
\sigma^{(j)} \sm \, \tau^{(j)}$ consists only of such roots and
(iii) $\tau \subseteq A$. As a consequence we have $\sigma \subseteq 
A$. Therefore $\sigma$ is an $(r-1)$-dimensional face of $\Delta$ 
containing $\tau$. This proves that $\Delta$ is pure of dimension 
$r-1$.

To prove that $\Delta$ is shellable we proceed by induction on $r$ and
the cardinality of $A$, the statement being trivial for $r \le 1$.
Let $S$ be the set of all pairs $(\alpha, i)$ such that $\alpha \in A$ has
color $i$, considered with the total ordering in which $(\alpha, i)$
precedes $(\alpha', i')$ if either $i<i'$ or $i=i'$ and $\alpha < \alpha'$
in the order of roots (\ref{total}). Theorem \ref{thm:Tz} implies that if
$(\alpha, i)$ is the smallest or the largest element of $S$ then
$\lk_\Delta (\alpha)$ coincides with the complex $\Delta^m_B (u)$, where
$u = w R(\alpha)$ in the former case and $u = R(\alpha) w$ in the latter,
and $B = A \cap \Phi^m_+ (u)$. Observe
that $B$ contains all roots in $\Phi^m_+ (u)$ of color $j$ and hence,
by the previous paragraph and the induction hypothesis, $\lk_\Delta
(\alpha)$ is pure of dimension $\ell (u) -1 = r-2$ and shellable. We now
distinguish two cases. If $A$ does not contain any vertices of color
other than $j$ then $\Delta$ is combinatorially isomorphic to $\Delta_+
(w)$ and hence shellable by Theorem \ref{thm:ABMW}. Otherwise the smallest
or the largest pair $(\alpha, i)$ in $S$ satisfies $i \neq j$. We already
know that $\lk_\Delta (\alpha)$ is pure $(r-2)$-dimensional and
shellable. By the same argument $\Delta \sm \, \alpha = \Delta^m_{A \sm
\, \{\alpha\}} (w)$ is pure $(r-1)$-dimensional and shellable, since $A
\sm \, \{\alpha\}$ contains all roots in $\Phi^m_+ (w)$ of color $j$.
Part (i) of Lemma \ref{lem0} implies that $\Delta$ is shellable as well.

\smallskip
\noindent
(ii) Let $\Delta^m_A (\Phi)$ denote the induced subcomplex of $\Delta^m
(\Phi)$ on the vertex set $A \subseteq \Phi^m_{\ge -1}$. Suppose that
$A$ contains all positive roots of color $j$ and set $\Delta =\Delta^m_A
(\Phi)$. We proceed by induction on the rank $n$ of $\Phi$
and the cardinality of $A$ and assume again, as we may, that $\Phi$ is
irreducible. If $A$ does not contain any negative simple root then the
result follows from the special case $w = \gamma$ of part (i). If
$\alpha = -\alpha_i \in A$ for some $\alpha_i \in \Pi$ then $\Delta
\sm \, \alpha = \Delta^m_{A \sm \, \{\alpha\}} (\Phi)$ is pure of
dimension $n-1$ and shellable by induction. Moreover, by Proposition
\ref{prop0} (ii) we have $\lk_\Delta (\alpha) = \Delta^m_B
(\Phi_\alpha)$, where $B = (\Phi_\alpha)_{\ge -1} \cap A$, so that
$\lk_\Delta (\alpha)$ is pure of dimension $n-2$ and shellable, again
by induction. Part (i) of Lemma \ref{lem0} implies that $\Delta$ is pure
of dimension $n-1$ and shellable.
\end{proof}

\smallskip
\noindent \emph{Proof of Theorem \ref{thm0}.}
(i) In view of Proposition \ref{prop0} (i) and Lemma \ref{lem0} (ii) we
may assume, once again, that $\Phi$ is irreducible.
Suppose that $B \subseteq \Phi^m_{\ge -1}$ has at most $m$ elements
and let $\Delta_A$ denote the induced subcomplex of $\Delta^m (\Phi)$
on the vertex set $A = \Phi^m_{\ge -1} \sm \, B$. If for some $1 \le i
\le m$ the set $B$ does not contain
any positive root of color $i$ then $\Delta_A$ is pure
$(n-1)$-dimensional and shellable by Proposition \ref{prop:pureshell}
(ii). Otherwise let $\beta \in B$ and, by part (b) of Proposition \ref{prop0}
(iii), choose $j$ so that $\rR_m^j (\beta) \in (-\Pi)$. Let $B' =
\rR_m^j (B)$ and $A' = \Phi^m_{\ge -1} \sm \, B'$. Part (a) of
Proposition \ref{prop0} (iii) implies that $\Delta_A$ is combinatorially
isomorphic to the induced subcomplex $\Delta_{A'}$ of $\Delta^m (\Phi)$
on the vertex set $A'$. Since $A'$ must contain all positive roots of
color $i$ for some $1 \le i \le m$, the result follows as in the first
part of the argument.

\smallskip
\noindent
(ii) Removing at most $m-1$ vertices from $\Delta^m_+ (\Phi)$ leaves all
positive roots of color $i$ in place for at least one index $1 \le i \le
m$. Hence the result follows from the special case $w = \gamma$ of
Proposition \ref{prop:pureshell} (i).
\qed

\begin{remark}
{\rm It is known \cite[Theorem 2.1 (a)]{Ba} that if a simplicial complex
$\Delta$ is $k$-Cohen-Macaulay then for any non-maximal face $\sigma \in
\Delta$ the complex $\lk_\Delta (\sigma)$ has at least $k$ vertices. On
the other hand any codimension one face of $\Delta^m (\Phi)$ is contained
in exactly $m+1$ facets \cite[Proposition 3.10]{FR2}. It follows that
$m+1$ is the largest integer $k$ for which $\Delta^m (\Phi)$ is
$k$-Cohen-Macaulay. A similar statement follows for $\Delta^m_+ (\Phi)$
since there exist codimension one faces of $\Delta^m_+ (\Phi)$ which are
contained in exactly $m$ facets of $\Delta^m_+ (\Phi)$ (any positive
facet of $\Delta^m (\Phi_\alpha)$, with $\alpha \in \Pi$, will be such
a codimension one face of $\Delta^m_+ (\Phi)$).}
\label{rem:best}
\end{remark}

\vspace{0.1 in}
The following corollary completes a proof of \cite[Conjecture 4.6]{AT}.
\begin{corollary}
The complex $\Delta^m_+ (\Phi)$ has the homotopy type of a wedge of
spheres of dimension $n-1$. If $\Phi$ is irreducible
then the number of these spheres is equal to
\[ N^+ (\Phi, m-1) = \prod_{i=1}^n \frac{e_i + (m-1)h - 1}{e_i + 1},
\]
where $h$ is the Coxeter number and $e_1, e_2,\dots,e_n$ are the
exponents of $\Phi$.
\label{cor:euler}
\end{corollary}
\begin{proof}
The first statement is part of Theorem \ref{thm0} (ii).
For the second statement it suffices to show that $\Delta^m_+ (\Phi)$ has
reduced Euler characteristic $\tilde{\chi} (\Delta^m_+ (\Phi)) =
(-1)^{n-1} N^+ (\Phi, m-1)$. This was shown in \cite{AT} for irreducible
root systems $\Phi$ of classical type with a proof that can be extended
to any type. We give the details here for the record.
To match the notation of \cite{AT} we write $\Phi = \Phi_I$, where $I$ is
an index set in bijection with $\Pi$, and denote by $\Phi_J$ the parabolic
root subsystem corresponding to $J
\subseteq I$. We also define $N^+ (\Phi_I, m)$ multiplicatively on the
irrecucible components of $\Phi_I$, if $\Phi_I$ is reducible. As was
mentioned in Remark 3 of \cite[Section 7]{AT}, Theorem 3.7 in \cite{FR2}
implies that equation (6) in \cite{AT} holds without any restriction
on $\Phi$. Multiplying both sides with $(-1)^{k-1}$ and summing over
$k$ we get
\[ (-1)^{|I|-1} \, \tilde{\chi} (\Delta^m (\Phi_I)) = \sum_{J \subseteq
I} \ (-1)^{|J|-1} \, \tilde{\chi} (\Delta^m_+ (\Phi_J)). \]
On the other hand
\[ (-1)^{|I|-1} \, \tilde{\chi} (\Delta^m (\Phi_I)) = \sum_{J \subseteq
I} \ N^+ (\Phi_J, m-1) \]
by Propositions 11.1 and 12.3 in \cite{FR2}. Inclusion-exclusion implies
the desired formula for $\tilde{\chi} (\Delta^m_+ (\Phi))$.
\end{proof}

\begin{remark}
{\rm Using results of \cite[Section 4]{ABMW} (in particular
\cite[Corollary 4.6]{ABMW}) it is possible to show that removing the
vertices of $\Delta_+ (w)$ in the order prescribed by (\ref{total})
gives a vertex decomposition of $\Delta_+ (w)$, in the sense of
\cite{PB}, for any $w \preceq \gamma$. It follows from the proofs of
our results in this section that $\Delta^m (\Phi)$ and $\Delta^m_+
(\Phi)$ are vertex decomposable as well. }
\label{rem:vert}
\end{remark}

\section{Generalized noncrossing partitions}
\label{gnc}

In this section we discuss some connections with the posets of 
generalized noncrossing partitions, defined and studied by D. 
Armstrong \cite{Ar2}. Let $\gamma \in W$ be any Coxeter element and
$\nc^{(m)}_W$ denote the set of $m$-tuples $(w_1, w_2,\dots,w_m)$ of
elements of $W$ satisfying $w_1 w_2 \cdots w_m \preceq \gamma$ and
$$\ell_T (w_1 w_2 \cdots w_m) = \ell_T (w_1) + \ell_T (w_2) + \cdots 
+ \ell_T (w_m).$$ We partially order this set by letting
\[ (u_1, u_2,\dots,u_m) \preceq (w_1, w_2,\dots,w_m) \ \ \text{if} 
\ \ u_i \preceq w_i \ \ \text{for} \ \ 1 \le i \le m \]
(this is the dual to the partial order defined in \cite{Ar2}). The
isomorphism type of $\nc^{(m)}_W$ does not depend on the choice of 
$\gamma$.
\begin{remark}
{\rm We have $(u_1, u_2,\dots,u_m) \in \nc^{(m)}_W$ whenever $(w_1,
w_2,\dots,w_m) \in \nc^{(m)}_W$ and $u_i \preceq w_i$ for all $1 
\le i \le m$.}
\label{rem:gnc}
\end{remark}
The poset $\nc^{(m)}_W$ has a unique minimal element $\hat{0}$ and 
is graded with rank function $\ell_T (w_1, w_2,\dots,w_m) := \ell_T 
(w_1 w_2 \cdots w_m)$. Now let $\gamma$ be as in Section \ref{pre4}. 
By Theorem \ref{thm:Tz} we have
\[ f(\sigma) = (w_{\sigma^{(m)}}, 
w_{\sigma^{(m-1)}},\dots,w_{\sigma^{(1)}}) \in \nc^{(m)}_W \]
for any face $\sigma$ of $\Delta^m_+ (\Phi)$. Recall that the 
\emph{upper truncation} $P^{\le k}$ of a graded poset $P$ is the 
induced subposet on the set of elements of rank at most $k$ (see 
\cite[Chapter 3]{Sta} for basic background on partially ordered 
sets). If $\Delta$ is a simplicial complex then a map $g: \Delta - 
\{\emptyset\} \to P$ is \emph{order preserving} if $g(\tau) \le 
g(\sigma)$ holds in $P$ whenever $\tau \subseteq \sigma$ are
nonempty faces of $\Delta$. In what follows, we say that a poset $P$ 
has a certain topological property if the geometric realization of 
the simplicial complex of chains in $P$ \cite[Section 9]{Bj} has
the same property.
\begin{theorem}
The map $$f: \Delta^m_+ (\Phi) - \{\emptyset\} \to \nc^{(m)}_W - 
\hat{0}$$ is well defined and order preserving and induces a 
homotopy equivalence between the skeleton $(\Delta^m_+ (\Phi))^{\le 
k-1}$ and $(\nc^{(m)}_W)^{\le k} - \hat{0}$ for all $1 \le k \le n$.

In particular, $\nc^{(m)}_W - \hat{0}$ is homotopy equivalent to
$\Delta^m_+ (\Phi)$.
\label{thm:f}
\end{theorem}
\begin{proof}
Let us write $\Delta = \Delta^m_+ (\Phi)$ and $\nc = \nc^{(m)}_W$, to
simplify notation. That the map $f$ is well defined and order preserving
follows from Theorem \ref{thm:Tz}. Moreover
we have $\ell_T (f(\sigma)) = |\sigma|$ for all $\sigma \in \Delta$ and
hence $f$ induces a well defined order preserving map from $\Delta^{\le
k-1}$ to $\nc^{\le k} - \hat{0}$ for all $1 \le k \le n$. To
complete the proof, by Quillen's Fiber Theorem \cite[Theorem 10.5 (i)]{Bj}
it suffices to show that $f^{-1} (\nc_{\le x})$ is contractible for all $x
\in \nc - \hat{0}$. Indeed, it follows directly from the definitions and
Remark \ref{rem:gnc} that if $x = (w_1, w_2,\dots,w_m)$ then $f^{-1}
(\nc_{\le x})$ is combinatorially isomorphic to the simplicial join of the
complexes $\Delta^m_+ (w_i)$ for $1 \le i \le m$, which is contractible by
Theorem \ref{thm:ABMW}.
\end{proof}

\begin{remark}
{\rm (i) It follows from Corollary \ref{cor:euler} and Theorem \ref{thm:f}
that $\nc^{(m)}_W - \hat{0}$ is homotopy equivalent to a wedge of
$(n-1)$-dimensional spheres. With some more work, and using the previous
statement as well as the main result of \cite{ABW}, one can show
that the poset $\nc^{(m)}_W \cup \hat{1}$, obtained from $\nc^{(m)}_W$ by
adjoining a maximum element $\hat{1}$, is homotopy Cohen-Macaulay (hence
Cohen-Macaulay). We omit
the details since a stronger statement, namely that $\nc^{(m)}_W \cup
\hat{1}$ is shellable, has been proved by D. Armstrong and H. Thomas
\cite[Theorem 3.7.2]{Ar2}.

\smallskip
\noindent
(ii) The proof of Theorem \ref{thm:f} shows that any order ideal $\iI$ in
$\nc^{(m)}_W - \hat{0}$ is homotopy equivalent to the subcomplex $f^{-1}
(\iI)$ of $\Delta^m_+ (\Phi)$.

\smallskip
\noindent
(iii) In the case $m=1$ and $k=n-1$, Theorem \ref{thm:f} gives a new proof
of the fact \cite[Corollary 4.4]{ABW} that the M\"obius number $\mu (\hat{0},
\gamma)$ of $\nc^{(1)}_W$ is equal, up to the sign $(-1)^n$, to the
number of facets of $\Delta_+ (\Phi)$ (positive clusters associated to
$\Phi$).  }
\label{rem:gnc2}
\end{remark}

\vspace{0.1 in} \noindent \emph{Acknowledgements}. The authors thank 
the anonymous referee for several helpful suggestions.

\end{document}